\documentclass[production,hyp]{nyj}
\recdate{October 4, 2006. Revised June 6, 2007}
\beginningpage{147}
\papnum{9}
\pubvol{13}
\pubyr{2007}

\renewcommand{\labelenumi}{(\arabic{enumi})}
\renewcommand{\theenumi}{\arabic{enumi}}
\makeatletter
\newcommand{\cmpl}[1]{%
   \sbox\z@{$#1$}%
   \dimen@=\wd\z@
   \advance \dimen@ -\strip@pt\fontdimen\@ne\textfont\@ne \ht\z@
   \setbox\tw@=\hb@xt@\dimen@{}%
   \ht\tw@=\ht\z@ \dp\tw@=\dp\z@
   \box\z@
   \llap{$\overline{\box\tw@}$}%
}
\newcommand{\cmpls}[1]{%
   \sbox\z@{$\scriptstyle #1$}%
   \dimen@=\wd\z@
   \advance \dimen@ -\strip@pt\fontdimen\@ne\textfont\@ne \ht\z@
   \setbox\tw@=\hb@xt@\dimen@{}%
   \ht\tw@=\ht\z@ \dp\tw@=\dp\z@
   \box\z@
   \llap{$\overline{\box\tw@}$}%
}
\makeatother

\makeatletter
\def\eqalignbot#1{%
 \null\,\vbox{\openup\jot\m@th
  \ialign{\strut\hfil$\displaystyle{##}$&$\displaystyle{{}##}$\hfil
      \crcr#1\crcr}}\,}
\makeatother

\newcommand{\twoshref}[1]{\href{#1two}{#1}}
\theoremstyle{plain}
\newtheorem{theorem}{Theorem}

\newtheorem{lemma}[theorem]{Lemma}

\theoremstyle{definition}

\newcommand{\R}{{\mathbb R}}
\newcommand{\Z}{{\mathbb Z}}
\newcommand{\F}{{\mathbb F}}

\title{Uniqueness of the $(22,891,1/4)$ spherical code}

\author[Henry Cohn]{\href{cohn}{Henry Cohn}}
\address{\hlabel{cohn}Microsoft Research, One Microsoft Way, Redmond, WA 98052-6399}
\email{\mailurl{cohn@microsoft.com}}

\author[Abhinav Kumar]{\href{kumar}{Abhinav Kumar}}
\address{\hlabel{kumar}Department of Mathematics, Harvard University, Cambridge, MA 02138}
\email{\mailurl{abhinav@math.harvard.edu}}
\curraddr{Microsoft Research, One Microsoft Way, Redmond, WA 98052-6399}
\email{\mailurl{abhinavk@microsoft.com}}

\thanks{Kumar was supported by a summer internship in the Theory Group at Microsoft Research and a Putnam Fellowship at Harvard University.}

\keywords{Spherical code, kissing configuration, spherical design, Leech lattice}

\subjclass{52C17, 05B40}

\begin{document}

\begin{abstract}
We use techniques of Bannai and Sloane to give a new proof that
there is a unique $(22,891,1/4)$ spherical code; this result is
implicit in a recent paper by Cuypers. We also correct a minor
error in the uniqueness proof given by Bannai and Sloane for the
$(23,4600,1/3)$ spherical code.
\end{abstract}

\maketitle
\tableofcontents
\section[Introduction]{\hlabel{sec1}Introduction}

An $(n,N,t)$ spherical code is a set of $N$ points on the unit
sphere $S^{n-1} \subset \R^n$ such that no two distinct points in
the set have inner product greater than $t$.  In other words, the
angles between them are all at least $\cos^{-1} t$.  The
fundamental problem is to maximize $N$ for a given value of $t$,
or equivalently to minimize $t$ given $N$.  Of course, for
specific values of $N$ and $t$, maximality of $N$ given $t$ is not
equivalent to minimality of $t$ given $N$, but complete solutions
of these problems for all parameter values would be equivalent.

Linear programming bounds are a powerful tool for proving upper
bounds on $N$ given $t$ (see \cite{DGS}, \cite{KL}, or Chapter~{9}
in \cite{CS}). In particular, they prove sharp bounds in a number
of important special cases, listed in \cite{L}.  Once a code has
been proved optimal, it is natural to ask whether it is unique up
to orthogonal transformations.  That is known in every case to
which linear programming bounds apply except for one infinite
family that is not always unique (see Appendix~A of \cite{CK} for
an overview).  However, one should not expect uniqueness to hold
in general for optimal spherical codes: for example, the $D_5$
kissing arrangement appears to be an optimal $40$-point spherical
code in $\R^5$, but there is at least one other $(5,40,1/2)$ code
(see \cite{Leech}).

One noteworthy case is the $(22,891,1/4)$ code.  A proof of
uniqueness is implicit in the recent paper \cite{Cu} by Cuypers,
but we are unaware of any explicit discussion of uniqueness in the
literature (by contrast, every other case has been explicitly
analyzed).  In this paper, we apply techniques from \cite{BS} to
give a new proof that it is unique.

This code arises naturally in the study of the Leech lattice in
$\R^{24}$ (see \cite{E} or \cite{CS} for background).  In the
sphere packing derived from the Leech lattice, each sphere is
tangent to $196560$ others.  The points of tangency form a
$(24,196560,1/2)$ code known as the kissing configuration of the
Leech lattice.  It can be viewed as a packing in $23$-dimensional
spherical geometry, whose kissing configuration is a
$(23,4600,1/3)$ code.  The $(22,891,1/4)$ code is obtained by
taking the kissing configuration once more; it is well defined
because the automorphism group of the Leech lattice acts distance
transitively on the $(24,196560,1/2)$ code.  All three of these
codes are optimal (in fact, universally optimal---see \cite{CK}),
although that is not known for the $(21,336,1/5)$ code that comes
next in the sequence.  The linear programming bounds are not sharp
for the $(21,336,1/5)$ code, and we make no conjecture as to
whether it is optimal. Its kissing configuration is a
$(20,170,1/6)$ code whose symmetry group does not even act
transitively: there are two orbits of points, one with $10$ points
(forming the midpoints of the edges of a regular $4$-dimensional
simplex) and one with $160$ points.  Because of the lack of
transitivity, this configuration has two different types of
kissing configurations, and it seems fruitless to continue
examining iterated kissing configurations. The $(20,170,1/6)$ code
is not universally optimal and probably not even optimal.

One can also construct the $(22,891,1/4)$ code using a
$6$-dimensional Hermitian space over $\F_4$.  Points in the
configuration correspond to $3$-dimensional totally isotropic
subspaces, with the inner product between two points ($-1/2$,
$1/4$, or $-1/8$) determined by the dimension of the intersection
of the corresponding subspaces ($2$, $1$, or $0$, respectively).
The graph with these subspaces as vertices and with edges between
pairs of subspaces with intersection dimension $2$ is the dual
polar graph associated with the group $\textup{PSU}(6,2)$ (see
Section~{9.4} in \cite{BCN}). In the paper \cite{Cu}, it is implicit
in the proof of Proposition~2.2 that a $(22,891,1/4)$ spherical
code must have the combinatorial structure of a $(2,4,20)$ regular
near hexagon, which is equivalent to this dual polar space
structure (see \cite{SY}).  Uniqueness then follows from the
classification of all polar spaces of rank at least $3$ by Tits in
\cite{T}.  By contrast, our proof makes use of entirely different
machinery.

The linear programming bounds not only prove bounds on spherical
codes, but also provide additional information about the codes
that achieve a given bound. When used with the auxiliary
polynomial $(x+1/2)^2(x+1/8)^2(x-1/4)$, they prove that every code
in $S^{21}$ with maximal inner product $1/4$ has size at most
$891$, and that equality is achieved iff all inner products
between distinct vectors are in $\{-1/2,-1/8,1/4\}$ and the code
is a spherical $5$-design.  Recall that a spherical $t$-design is
a finite subset of the sphere $S^{n-1} \subset \R^n$ such that for
every polynomial function $p \colon \R^{n} \to \R$ of total degree
at most $t$, the average of $p$ over the design equals its average
over the entire sphere.

The techniques we use to prove uniqueness were developed by Bannai
and Sloane in \cite{BS}, and we follow their approach quite
closely. (Note that their paper is reprinted as Chapter~{14} of
\cite{CS}.) They proved uniqueness for the $(24,196560,1/2)$ and
$(23,4600,1/3)$ codes, as well as analogous codes derived from the
$E_8$ root lattice.  Here we correct a minor error in their proof
for the $(23,4600,1/3)$ code.  They construct a lattice $L$ and
conclude their proof by saying ``and hence $L$ must be the Leech
lattice,'' but in fact it is not the Leech lattice (it is a
sublattice of index~$2$).  At the end of this paper we explain the
problem and how to correct it.

One small difference between this paper and \cite{BS} is that the
$(22,891,1/4)$ code is not a tight spherical design, whereas all
the designs dealt with in \cite{BS} are tight.  (A tight spherical
$(2e+1)$-design in $\R^n$ is one with $2\binom{n+e-1}{n-1}$
points, which by Theorem~5.12 of \cite{DGS} is a lower bound for
the number of points.) However, no fundamental changes in the
techniques are needed. The only important difference is that we
cannot conclude that the $(22,891,1/4)$ code is the only
$891$-point spherical $5$-design in $\R^{22}$, as we could if it
were tight.

\section[Uniqueness of the $(22,891,1/4)$ code]{\hlabel{sec2}Uniqueness of the $(22,891,1/4)$ code}

\begin{theorem} \bothlabel{theorem:891}
There is a unique $(22,891,1/4)$ spherical code, up to orthogonal
transformations of $\R^{22}$.
\end{theorem}

Let $\mathcal{C}$ be such a code.  We begin with the observation
that by the sharpness of the linear programming bounds, $-1/2$,
$-1/8$, and $1/4$ are the only possible inner products that can
occur between distinct points in $\mathcal{C}$.  Let
$u_1,\dots,u_{891}$ be the points in $\mathcal{C}$, and let
\begin{align*}
U_i &= (1, 1/\sqrt{3}, \sqrt{8/3}\, u_i), \\
V_0 &= (2,0,\dots, 0), \textup{ and} \\
V_1 &= (1,\sqrt{3},0, \dots, 0)
\end{align*}
be vectors in $\R^{24}$.  The slightly nonstandard notation $(1,
1/\sqrt{3}, \sqrt{8/3}\, u_i)$ of course means the concatenation
of the vectors $(1,1/\sqrt{3})$ and $\sqrt{8/3}\, u_i$.

It is easy to check that all these vectors have norm $4$ and the
inner product between any two of them is an integer; specifically,
$\langle U_i,U_j \rangle$ is $4$, $2$, $1$, or $0$ according as
$\langle u_i, u_j \rangle$ is $1$, $1/4$, $-1/8$, or $-1/2$,
respectively. Let $L$ be the lattice spanned by
$U_1,\dots,U_{891}$, $V_0$, and $V_1$. It follows that $L$ is an
even lattice (i.e., all vectors have even norms). We will show
that $L$ is uniquely determined, up to orthogonal transformations
of $\R^{24}$ that fix $V_0$ and $V_1$, as is
$\{U_1,\dots,U_{891}\}$.

In what follows, vectors in $\R^{24}$ are generally denoted by
uppercase letters and vectors in $\R^{22}$ by lowercase letters.
One exception is the standard basis $e_1,\dots,e_{24}$ of
$\R^{24}$.

\begin{lemma} \bothlabel{lemma:norm4}
The minimal norm $\langle V,V \rangle$ for $V \in L \setminus
\{0\}$ is $4$.
\end{lemma}

\begin{proof}
Suppose there exists $V \in L$ with $\langle V, V\rangle = 2$.
Then $\langle V,W \rangle \in \{0,\pm 1, \pm 2\}$ for all $W \in
\{V_0,V_1,U_1,\dots,U_{891}\}$, because $\langle V,W \rangle \in
\mathbb{Z}$ and $|\langle V,W \rangle| \leq |V| |W| = 2\sqrt{2}$.

Now let $V = (x,y/\sqrt{3},\sqrt{8/3}\,u)$ with $u \in \R^{22}$
and $x,y \in \R$. We note that $x$ and $y$ must be integers of the
same parity, from the description of the generators of the lattice
$L$. Also, we must have $x^2 + y^2/3 \leq 2$, by the condition on
the norm of $V$. This implies that $(x,y) \in \{(0,0), (0, \pm 2),
(\pm 1, \pm 1) \}$.  We can furthermore assume that $(x,y) \in
\{(0,0), (0, 2), (1, \pm 1) \}$, because otherwise we replace $V$
with $-V$.  If $(x,y) = (0,2)$, then $\langle V,V_1 \rangle = 2$
and thus $|V_1-V|^2=2$, so we can replace $V$ with $V_1-V$ and
$(x,y)$ with $(1,1)$. If $(x,y) = (1,-1)$, then we can replace $V$
with $V_0-V$ and $(x,y)$ with $(1,1)$.  We can therefore assume
that $(x,y)$ is $(0,0)$ or $(1,1)$.

If $(x,y)=(1,1)$, then we claim that there exists an $i$ such that
$\langle V, U_i \rangle = 2$, in which case replacing $V$ with
$V-U_i$ reduces to the case of $(x,y)=(0,0)$.  To prove the
existence of such an $i$, consider the point $u \in \R^{22}$,
which has $|u|=1/2$. For each $i$, if $\langle V, U_i \rangle \in
\{-2,-1,0,1\}$, then $\langle u,u_i \rangle \in
\{-5/4,-7/8,-1/2,-1/8\}$.  If that were always the case, then the
set $\{2u,u_1,\dots,u_{891}\}$ would be a $(22,892,1/4)$ spherical
code, which is impossible.

We are left with only one case, namely that $(x,y) = (0,0)$. Then
$|u| = \sqrt{3}/{2}$. The inner products $\langle u, u_i \rangle$
must lie in the set $\{0 ,\pm 3/8, \pm 3/4\}$, corresponding to
the restriction that $\langle V, U_i \rangle \in \{0, \pm 1, \pm
2\}$. Let $N_{0}, N_{3/8}, N_{-3/8}, N_{3/4}, N_{-3/4}$ be the
numbers of vectors $u_i$ that have inner products
$0,3/8,-3/8,3/4,-3/4$, respectively, with $u$. Now from the fact
that $\mathcal{C}$ is a $5$-design (which we obtain from the
linear programming bounds), we observe that for every polynomial
$p(x)$ of degree at most $5$,
$$
\frac{\sum_{\alpha \in \{0,3/8,-3/8,3/4,-3/4\}} N_\alpha
p(\alpha)}{891} = \int_{S^{21}} p\big(\langle z, u \rangle\big)
\,d\mu(z),
$$
where the surface measure $\mu$ on $S^{21}$ has been normalized to
have total volume $1$.

The right side does not depend on the direction of $u$, only on
its magnitude, and it is easily evaluated when $p(x)=x^i$: for $i$
odd it vanishes, and for $i$ even it equals
$$
|u|^i \frac{i!(22/2-1)!}{(i/2+22/2-1)!(i/2)!2^i} =
\frac{i!\,10!}{(10+i/2)!(i/2)!}
\left(\frac{\sqrt{3}}{4}\right)^{i}.
$$

We write down five equations corresponding to the monomials
$p(x)=1$, $x$, $x^2$, $ x^3$, and $x^4$ and solve the resulting
system of equations to get
$$
(N_0, N_{3/8}, N_{-3/8}, N_{3/4}, N_{-3/4}) = (657, 120, 120, -3,
-3).
$$
The negative numbers give us the contradiction.
\end{proof}

As an immediate corollary we observe that the (integral) inner
product between two minimal vectors of $L$ cannot be $\pm 3$ and
so must lie in $\{0,\pm 1, \pm 2, \pm 4\}$: if $\langle U,V
\rangle = 3$ with $U$ and $V$ minimal vectors, then $|U-V|^2 =
|U|^2 + |V|^2 - 2\langle U,V \rangle = 2$, contradicting
Lemma~\ref{lemma:norm4}.

\begin{table}\hlabel{table:intnums}
\begin{align*}
P_1(1/4,1/4) &= 336 & P_1(-1/8,-1/8) &= 512 & P_1(-1/2,-1/2) &=  42\\[3pt]
P_{1/4}(1/4,1/4) &= 170 & P_{1/4}(-1/8,-1/8) &= 320 & P_{1/4}(-1/2,-1/2) &= 5\\
P_{1/4}(1/4,-1/8) &= 160 & P_{1/4}(1/4,-1/2) &= 5 & P_{1/4}(-1/8,-1/2) &= 32\\[3pt]
P_{-1/8}(1/4,1/4) &= 105 & P_{-1/8}(-1/8,-1/8) &= 280 & P_{-1/8}(-1/2,-1/2) &= 0\\
P_{-1/8}(1/4,-1/8) &= 210 & P_{-1/8}(1/4,-1/2) &= 21 & P_{-1/8}(-1/8,-1/2) &= 21\\[3pt]
P_{-1/2}(1/4,1/4) &= 40 & P_{-1/2}(-1/8,-1/8) &= 256 & P_{-1/2}(-1/2,-1/2) &= 1\\
P_{-1/2}(1/4,-1/8) &= 256 & P_{-1/2}(1/4,-1/2) &= 40 & P_{-1/2}(-1/8,-1/2) &= 0\\
\end{align*}
\caption{Intersection numbers for a $(22,891,1/4)$ code.}
\label{table:intnums}
\end{table}

It follows from Theorem~7.4 in \cite{DGS} that because
$\mathcal{C}$ is a $5$-design in which $3$ inner products other
than $1$ occur and $5 \ge 2\cdot3 - 2$, the points in
$\mathcal{C}$ form a $3$-class association scheme when pairs of
points are grouped according to their inner products.  In other
words, given $\alpha,\beta,\gamma \in \{-1/2,-1/4,1/8,1\}$, there
is a number $P_\gamma(\alpha,\beta)$ such that whenever $\langle
u_i,u_j \rangle = \gamma$, there are exactly
$P_\gamma(\alpha,\beta)$ points $u_k$ such that $\langle u_i,u_k
\rangle = \alpha$ and $\langle u_j,u_k \rangle = \beta$.  These
numbers are called intersection numbers and are determined in the
proof of the theorem in \cite{DGS}.  We have tabulated them in
Table~\ref{table:intnums} (note that $P_\gamma(\alpha,\beta) =
P_\gamma(\beta,\alpha)$, $P_\gamma(\alpha,1)$ is the Kronecker
delta function $\delta_{\alpha,\gamma}$, and $P_1(\alpha,\beta) =
0$ unless $\alpha=\beta$).

\begin{lemma} \bothlabel{lemma:Dn}
The lattice $L$ contains a sublattice isometric to
$\sqrt{2}D_{24}$ and containing $V_0$ and $V_1$.
\end{lemma}

Recall that the minimal norm in $D_n$ is $2$, so it is $4$ in
$\sqrt{2}D_n$.

\begin{proof}
We prove by induction on $n$ that there exist minimal vectors
$G_1,\dots,G_n \in L$ such that $\langle G_1, G_2 \rangle = 0$,
$\langle G_1, G_3 \rangle = -2$, $\langle G_i,G_{i+1}\rangle = -2$
for $2 \le i \le n-1$, and all other inner products vanish.  In
other words, for $3 \le k \le n$, the vectors $G_1,\dots,G_k$ span
a copy of $\sqrt{2}D_k$, as one can see from the Dynkin diagram of
$D_k$:
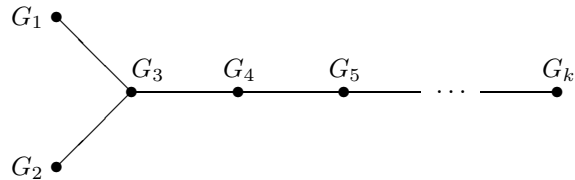
\begin{center}
\setlength{\unitlength}{1cm}
\begin{picture}(7.55684,2.4)(-0.7,-0.2)
\put(0,0){\circle*{0.15}} \put(1,1){\circle*{0.15}}
\put(0,2){\circle*{0.15}} \put(2.41421,1){\circle*{0.15}}
\put(3.81842,1){\circle*{0.15}} \put(6.65684,1){\circle*{0.15}}
\put(1,1){\line(-1,1){1}} \put(1,1){\line(-1,-1){1}}
\put(1,1){\line(1,0){1.4142135}}
\put(2.41421,1){\line(1,0){1.4142135}}
\put(3.81842,1){\line(1,0){1.0142135}}
\put(6.65684,1){\line(-1,0){1.0142135}}
\put(5.24263,1){\raisebox{-0.0175cm}{\hskip -0.2cm $\dots$}}
\put(-0.2,2){\raisebox{-3pt}{\hskip -0.4cm $G_1$}}
\put(-0.2,0){\raisebox{-3pt}{\hskip -0.4cm $G_2$}}
\put(1,1.2){$G_3$} \put(2.41421,1.2){\hskip -0.2cm $G_4$}
\put(3.81842,1.2){\hskip -0.2cm $G_5$} \put(6.65684,1.2){\hskip
-0.2cm $G_k$}
\end{picture}
\end{center}
In what follows we refer to this copy as $\sqrt{2}D_k$, and we
write $G_1 = -\sqrt{2}(E_1+E_2)$ and $G_i =
\sqrt{2}(E_{i-1}-E_{i})$ for $i \ge 2$, so $E_1,\dots,E_n$ is an
orthonormal basis of the ambient space $\R D_n = \R \otimes_\Z
\sqrt{2}D_n$ of $\sqrt{2} D_n$.

We will furthermore choose this sublattice to contain $V_0$ and
$V_1$ when $n \ge 5$.

For $n=3$, the existence of such vectors follows immediately from
the fact that the intersection numbers $P_1(-1/2,-1/2)=42$ and
$P_{-1/2}(1/4,1/4)=40$ are positive.  Choose $G_1=U_i$ for any
$i$.  Then among $U_1,\dots,U_{891}$ there are $42$ choices for
$G_2$, and $40$ choices among $-U_1,\dots,-U_{891}$ for $G_3$
given $G_2$.

Now suppose the assertion holds up to dimension $n$, with $3 \le n
< 24$.  As a first step we show that there are at least $43$
minimal vectors $W$ in $L$ such that $\langle G_i, W \rangle = 2$
for $i \in \{1,2\}$, whereas in $\sqrt{2}D_n$ there are only $2n-4
\le 42$, namely $G_1+G_2 + \dots +G_k$ and $-G_3-G_4-\dots-G_k$
with $3 \le k \le n$.  (Checking this assertion for $\sqrt{2}D_n$
is a straightforward exercise in manipulating coordinates.)  The
next lattice $\sqrt{2}D_{n+1}$ will be spanned by $\sqrt{2}D_n$
and such a vector $W$.

To construct these vectors $W$ we work as follows.  Renumbering
the vectors if necessary, we can assume that $U_1,\dots,U_{40}$
satisfy $\langle G_i, U_j \rangle = 2$ for $i \in \{1,2\}$ and $1
\le j \le 40$ (because $P_{-1/2}(1/4,1/4)=40$).  The vectors $V_0$
and $V_1$ also satisfy $\langle G_i, V_j \rangle = 2$ for $i \in
\{1,2\}$ and $j \in \{0,1\}$. We must still find one more choice
of $W$. To do so, note that $P_{-1/2}(-1/2,-1/2)=1$.  Hence there
is a unique vector $U_\ell$ such that $\langle U_\ell, G_1\rangle
= \langle U_\ell, G_2\rangle = 0$.  The vector $V_2 = V_0 -
U_\ell$ is another choice for $W$ (we could also choose
$V_1-U_\ell$, but we will not require that many possibilities).

The $43$ vectors $U_1,\dots,U_{40},V_0,V_1,V_2$ are all distinct:
the only possible danger is if $V_2$ equals one of the other
vectors.  Because $V_2=V_0-U_\ell$, clearly $V_2 \ne V_0$, and
$V_2 \ne V_1$ follows from looking at the second coordinate in the
definitions of $V_0,V_1,U_i$.  Similarly, $V_2 = U_i$ is
impossible because comparing second coordinates shows that $V_0
\ne U_i+U_\ell$.

Thus, there are at least $43$ minimal vectors $W$ satisfying
$\langle G_k, W \rangle = 2$ for $k = 1,2$, whereas in
$\sqrt{2}D_n$ there are at most $42$.  Choose a minimal vector $W$
with this property such that $W \notin \sqrt{2}D_n$, and in
particular choose $W=V_0$ or $W=V_1$ if possible.  (That will
ensure that $V_0,V_1 \in \sqrt{2}D_n$ if $n \ge 5$.)

This vector $W$ cannot be in $\R D_n$: if $W = \sum_{i=1}^n c_i
E_i$, then $\langle G_k, W \rangle = 2$ for $k \in \{1,2\}$
implies $c_1 =0$ and $c_2=-\sqrt{2}$. For $3 \leq i \leq n$,
$\sqrt{2}(E_2\pm E_i)$ is a minimal vector in $\sqrt{2}D_n \subset
L$, and therefore $\langle W, \sqrt{2} (E_2 \pm E_i)\rangle \in
\{0, \pm 1, \pm 2\}$.  (The inner product cannot be $\pm 4$
because $\sqrt{2}(E_2\pm E_i)\in \sqrt{2}D_n$ but
$W\not\in\sqrt{2}D_n$.) Because $\langle W, \sqrt{2} (E_2 \pm
E_i)\rangle = -2 \pm \sqrt{2}c_i$, it follows that $c_3 = c_4 =
\dots = c_n = 0$, which contradicts $\langle W,W \rangle = 4$.

Choose $E_{n+1}$ so that $\{E_1,\dots, E_{n+1}\}$ is an
orthonormal basis for $\R D_n \oplus \R W$, and let $W = c_1 E_1 +
\dots + c_{n+1} E_{n+1}$. Then the same calculation gives $c_1 =
0$, $c_2 = -\sqrt{2}$, $c_3 = \dots = c_n = 0$, and $c_{n+1} = \pm
\sqrt{2}$.  Thus, $\sqrt{2}D_n$ and $W$ span a copy of
$\sqrt{2}D_{n+1}$ contained in $L$.
\end{proof}

It will be convenient in the rest of the proof of
Theorem~\ref{theorem:891} to change coordinates to agree with the
standard coordinates for the Leech lattice (see
\cite[p.~131]{CS}).  To do so, choose coordinates so that $L$
contains the usual lattice $\sqrt{2}D_{24}$ (i.e., $\sqrt{2}$
times all the integral vectors with even coordinate sum), with
$V_0,V_1 \in \sqrt{2}D_{24}$. We can furthermore assume that $V_0
= (4,4,0,\dots, 0)/\sqrt{8}$ and $V_1 =
(4,0,4,0,0,\dots,0)/\sqrt{8}$, because the automorphism group of
$\sqrt{2}D_{24}$ acts transitively on pairs of minimal vectors
with inner product $2$.  (We write the vectors in this way, with
$4/\sqrt{8}$ instead of $\sqrt{2}$, because it will prove helpful
in dealing with the Leech lattice.)

In these new coordinates, all inner products are of course
preserved, but the coordinates for $U_1,\dots,U_{891}$ are no
longer the same as those we previously used.  Let
$e_1,\dots,e_{24}$ denote the standard basis of $\R^{24}$ with
respect to the new coordinates.  From this point on, all uses of
coordinates refer to the new coordinates.

We wish to show that the vectors $U_1,\dots,U_{891}$ are uniquely
determined, up to orthogonal transformations of $\R^{24}$ fixing
$V_0$ and $V_1$, which include of course permutations and sign
changes of the last $21$ coordinates. Let $W = (w_1,\dots,
w_{24})/\sqrt{8}$ be one of the $U_i$'s. Then
$$
\sum_{i=1}^{24} w_i^2 = 8|W|^2 = 32,
$$
$$ (w_i \pm w_j)/2 = \big\langle W, \sqrt{2}(e_i \pm e_j) \big\rangle \in
\{0, \pm 1, \pm 2, \pm 4\}
$$
for $i \ne j$,
$$
(w_1+w_2)/2 = \langle W, V_0 \rangle = 2,
$$
and
$$
(w_1+w_3)/2 = \langle W, V_1 \rangle = 2.
$$

{}From the above conditions we see that each $w_i$ is an integer
(because $(w_i+w_j)/2$ and $(w_i-w_j)/2$ are), and that they are
all at most $4$ in absolute value and of the same parity. A little
more work shows that the only possibilities are
\begin{equation}\bothlabel{eq:cases}
\sqrt{8}W = \begin{cases}\hlabel{I} 4(e_1\pm e_j) &\textup{with $j \ge 4$}, \\
\hlabel{II}
4(e_2+e_3), \\\hlabel{III}
2(e_1+e_2+e_3) +2 \sum_{k=1}^5 \pm e_{j_k} &\textup{with $3 < j_1
< j_2 < \dots < j_5$, or} \\\hlabel{IV}
3e_1 + e_2 + e_3 + \sum_{j=4}^{24} \pm e_j.
\end{cases}
\end{equation}
To prove this, note first that $w_1 \ge 0$ since $(w_1+w_2)/2=2$
and $w_2 \le 4$.  If $w_1=0$, then $w_2=w_3=4$, and $w_i=0$ for
$i>3$ because $|W|^2=4$.  If $w_1=1$, then $w_2=w_3=3$ and hence
$(w_2+w_3)/2 = 3$, which is impossible. If $w_1=2$, then
$w_2=w_3=2$; the constraint that $(w_1\pm w_i)/2 \in \{0, \pm 1,
\pm 2, \pm 4\}$ rules out $w_i = \pm 4$, so all remaining
coordinates are in $\{0, \pm 2\}$, and there must be five more
$\pm 2$'s because $|W|^2=4$. If $w_1=3$, then $w_2=w_3=1$, and
$(w_1\pm w_i)/2 \in \{0, \pm 1, \pm 2, \pm 4\}$ rules out $w_i =
\pm 3$ for $i>1$, so all remaining coordinates must be $\pm 1$.
Finally, if $w_1=4$, then $w_2=w_3=0$, and $(w_1\pm w_i)/2 \in
\{0, \pm 1, \pm 2, \pm 4\}$ implies $w_i \in \{0, \pm 4\}$;
exactly one more coordinate must be $\pm 4$ because $|W|^2=4$.

Call the cases enumerated in Equation \eqref{eq:cases} above
Case~\shref{I}, Case~\shref{II}, Case~\shref{III}, and Case~\shref{IV}, respectively.

By abuse of notation, view $\{0,1\}^{21}$ as being contained in
$\Z^{21} = \sum_{i=4}^{24} \Z e_i$. We define a code $\mathcal{D}
\subset \{0,1\}^{21}$ by stipulating that $c \in \mathcal {D}$ iff
$(2(e_1+e_2+e_3) + 2c+4z)/\sqrt{8}$ is one of the $U_i$'s for some
$z \in \Z^{21}$. This corresponds to Case~\shref{III} above. The codewords
in $\mathcal{D}$ have weight $5$, and the minimum distance between
codewords is at least $8$, since the minimum distance between
vectors of the lattice $L$ is $2$.  (If $(2(e_1+e_2+e_3) +
2c_1+4z_1)/\sqrt{8}$ and $(2(e_1+e_2+e_3) + 2c_2+4z_2)/\sqrt{8}$
are both as above, then $(2(c_1-c_2)+4(z_1-z_2))/\sqrt{8} \in L$.
One can add an element of $\sqrt{2}D_{24}$ to cancel all of
$4(z_1-z_2)/\sqrt{8}$ except for one coordinate, and another to
cancel the remaining coordinate at the cost of changing the sign
of one of the $\pm 2$'s occurring in $2(c_1-c_2)/\sqrt{8}$.  Then
if the distance between the codewords $c_1$ and $c_2$ in
$\mathcal{D}$ is less than $8$, the resulting vector in $L$ has
length less than $2$.)

It follows from the linear programming bounds for constant-weight
binary codes (see \cite[p.~545]{MS}) that the largest such code
has size $21$. In particular, it is a projective plane over $\F_4$
(the points are coordinates and the lines are the supports of the
codewords), or equivalently an $S(2,5,21)$ Steiner system, and it
is thus unique up to permutations of the coordinates (Satz~1 in
\cite{W}). Also, for each codeword of $\mathcal{D}$, we can only
use at most half of the possible sign assignments in the $\pm 2$'s
in Case~\shref{III}, since otherwise we would get two elements of $L$ that
agree except for one sign and are thus at distance
$(2-(-2))/\sqrt{8} = \sqrt{2}$, which is again a contradiction.
This gives a total of at most $2^4 \cdot 21 = 336$ possible
minimal vectors for Case~\shref{III}.

Similarly, for Case~\shref{IV}, define a code $\mathcal{E} \subset
\{0,1\}^{21}$ so that $c \in \mathcal{E}$ iff
$$
\left( 3e_1 + e_2 + e_3 + 2c -\sum_{i=4}^{24} e_i\right)/\sqrt{8}
$$
is one of the $U_i$'s. We note as before that codewords have
distance at least $8$ from each other, and also at most $16$
(otherwise two $U_i$'s would be too far apart). The largest such
code has $512$ codewords, as is easily proved using linear
programming bounds (see Theorem~20 of Chapter~{17} in
\cite[p.~542]{MS}), if one takes into account both the minimal and
the maximal distance.  This is more subtle than it might at first
appear, because the linear programming bounds are not in fact
sharp if one uses only the minimal distance. We conclude that
there are at most $512$ vectors in Case~\shref{IV}.

In all, the number of possible $U_i$'s is at most
$2\cdot21+1+336+512 = 891$. On the other hand, we already know
that there are $891$ of them. This forces the codes $\mathcal{D}$
and $\mathcal{E}$ to have the greatest possible size. In
particular, $\mathcal{D}$ is uniquely determined, up to
permutation of coordinates.  These coordinate permutations are
orthogonal transformations of $\R^{24}$ that fix $V_0$ and $V_1$
and preserve $\sqrt{2} D_{24}$. To complete the proof of
uniqueness, it will be enough to show that after performing
further such transformations that preserve the code $\mathcal{D}$,
we can specify all the vectors of Cases~\shref{III} and \shref{IV} exactly.

Let $W_0 = \left(3e_1 + e_2 + e_3 +2c_0 -\sum_{i=4}^{24}
e_i\right)/\sqrt{8}$ be a fixed vector from Case~\shref{IV}. Let
$i_1,\dots,i_r$ be the places (between $4$ and $24$) where $c_0$
has a $1$. Then let $\phi$ be the composition of reflections in
the corresponding hyperplanes (i.e., change the signs of those
coordinates). Applying $\phi$ clearly fixes $V_0$ and $V_1$ and it
takes the vector $W_0$ to $\left(3e_1 + e_2 + e_3 -\sum_{i=4}^{24}
e_i\right)/\sqrt{8}$. It also preserves the code $\mathcal{D}$.
Thus, we can assume that $W_0 = \left(3e_1 + e_2 + e_3
-\sum_{i=4}^{24} e_i\right)/\sqrt{8}$. Now we try to determine the
precise form of the vectors of Case~\shref{III}. We know that they have
$\pm 2$ entries in the positions of the code $\mathcal{D}$; the
only question is if we can pin down the positions of the signs.
Let $d \in \mathcal{D}$ be a codeword, and let $V$ be any vector
in Case~\shref{III} with $\pm 2$'s at the positions specified by the
codeword $d$.  Suppose $r$ of these are $-2$'s and $5-r$ are
$2$'s. Then taking the inner product with $W_0$, we get
$$
\langle W_0, V \rangle  = \frac{1}{8} (6+2+2 + 2r - 2(5-r)) =
\frac{4r}{8}.
$$
Since this inner product is an integer, we deduce that $r$ is
even. For each codeword $d$, this gives $2^5/2 = 2^4 = 16$
possible vectors. Thus the maximum number of allowed vectors is
$21 \cdot 16 = 336$, which we already know is the number of
vectors from Case~\shref{III}. Therefore equality holds, and we have
specified all the vectors of Case~\shref{III}. Namely, they are all
vectors of the form
$$
2e_1+2e_2+2e_3+ 2\sum_{j \textup{ such that } d_j = 1} \pm e_j,
$$
where an even number of minus signs are used and $d$ ranges over
all codewords in $\mathcal{D}$.

Now we claim that the lattice $L$ is generated by $\sqrt{2}
D_{24}$, the vectors in Case~\shref{III}, and $W_0$, which implies that
the vectors in Case~\shref{IV} are uniquely determined.  (Recall that they
are the only remaining vectors in $L$ that satisfy the constraints
enumerated in the paragraph before Equation \eqref{eq:cases}.)  To
show that $L$ is generated, it suffices to show that the vectors
in Case~\shref{IV} are, because all other generators are already included.

For this, let $W = \left(3e_1+ e_2 +e_3  + 2c - \sum_{i=4}^{24}
e_i \right)/\sqrt{8}$ be any vector in Case~\shref{IV}. Then $W-W_0 =
2c/\sqrt{8}$ is in the lattice $L$ and it is enough to show that
it is in the span of the above generators excluding $W_0$.
Equivalently, we must show that $c$ is in the span of $2(e_i \pm
e_j)$ and $e_1+e_2+e_3+d$ with $d \in \mathcal{D}$.  Because
$2c/\sqrt{8} \in L$ and $L$ is even, the weight of $c$ must be a
multiple of $4$.  Therefore, what we need to show is that in
$\F_2^{24}$, the codeword $000c$ is in the span of the codewords
$111d$ for $d \in \mathcal{D}$, where of course $000c$ denotes the
concatenation of $(0,0,0)$ with $c$. (When we work modulo $2$, the
vectors $2(e_i \pm e_j)$ vanish. Fortunately, that is not a
problem, because $000c$ and $111d$ all have weights divisible by
$4$.  It follows from congruence modulo $2$ that the difference in
$\Z^{24}$ between $000c$ and a sum of vectors of the form $111d$
is not only in the span of the vectors $2e_i$ but in fact in the
span of $2(e_i \pm e_j)$.)

Conversely, any vector of the form $000c$ that is in the span of
the codewords $111d$ for $d \in \mathcal{D}$ will correspond to a
vector in Case~\shref{IV}.

Of course one must take the sum of an even number of words $111d$
to arrive at a word of the form $000c$.  It is easily checked that
the code $\mathcal{D}$ spans a $10$-dimensional subspace of
$\F_2^{21}$ (simply check that the incidence matrix of the
projective plane over $\F_4$ has rank $10$ over $\F_2$; this is
easily checked directly or by using a general formula that is
implicit in \cite{GM} and explicit in \cite{MM} and \cite{S}).
Hence the codewords of the form $111d$ with $d \in \mathcal{D}$
span $512$ words of the form $000c$. Converting back to vectors,
these give us $512$ vectors of the form $W-W_0$ with $W$ in
Case~\shref{IV}. However, we know that the total number of $W$'s in
Case~\shref{IV} is $512$. Therefore all of them must come from this
construction. In other words, this shows that $000c$ is always in
the span of $111d$ with $d \in \mathcal{D}$, which is what we
wanted to prove.

This concludes the proof of Theorem~\ref{theorem:891}.

\section[Uniqueness of the $(23,4600,1/3)$ code]{\hlabel{sec3}Uniqueness of the $(23,4600,1/3)$ code}

Now we add a correction to the proof in \cite{BS} that there is a
unique code of size $4600$ and maximum inner product $1/3$ in
$S^{22}$, up to orthogonal transformations of $\R^{23}$. As
mentioned above, this code is derived from the Leech lattice by
taking the kissing arrangement twice.

\begin{theorem} \bothlabel{theorem:4600}
There is a unique $(23,4600,1/3)$ spherical code, up to orthogonal
transformations of $\R^{23}$.
\end{theorem}

\begin{proof}
Let $u_1,\dots,u_{4600}$ be the points in the code, and set
$$
V_0 = (2,0,\dots,0) \in \R^{24}
$$
and
$$
U_i = (1, \sqrt{3} u_i).
$$
Let $L$ be the lattice spanned by $V_0$ and $U_1,\dots,U_{4600}$.

The analogues of Lemmas~\ref{lemma:norm4} and~\ref{lemma:Dn} go
through as before. However, it is then stated in \cite{BS} that
$L$ is the Leech lattice, which is not correct (for by
construction, every element of $L$ has even inner product with
$V_0$, which is not true for every vector in the kissing
configuration of the Leech lattice). However, one can take the
path that we have described above. Briefly, we have the following
setup:

Choose new coordinates so that $L$ contains the usual lattice
$\sqrt{2}D_{24}$ and $V_0 = (4e_1+4e_2)/\sqrt{8}$. The vectors in
$L$ that could possibly have inner product $2$ with $V_0$ are of
the form $(w_1,\dots,w_{24})/\sqrt{8}$ with
$$
(w_1,w_2,\dots,w_{24}) = \begin{cases}\hlabel{Itwo}
 4(e_1+e_j) & \textup{with $j \geq 3$,} \\\hlabel{IItwo}
4(e_2+e_j) & \textup{with $j \geq 3$,}\\\hlabel{IIItwo}
2(e_1+e_2) +2 \sum_{k=1}^6 \pm e_{j_k} & \textup{with $2 < j_1 <
j_2 < \dots < j_6$,} \\\hlabel{IVtwo}
 3e_1 + e_2  + \sum_{j=3}^{24} \pm e_j, & \textup{or} \\\hlabel{Vtwo}
e_1 + 3e_2  + \sum_{j=3}^{24} \pm e_j.
\end{cases}
$$
Call these cases Case~\twoshref{I} through Case~\twoshref{V}.

Once again we enumerate the possibilities: Cases~\twoshref{I} and~\twoshref{II} lead to
$22$ vectors each. Case~\twoshref{III} leads to a $(22,8,6)$ code, which has
at most $77$ elements by the linear programming bounds for
constant-weight codes. Therefore there are at most $2^5 \cdot 77 =
2464$ vectors from Case~\twoshref{III} (as in the previous case only half of
the possible sign patterns can occur). Finally, Cases~\twoshref{IV} and~\twoshref{V}
both lead to $(22,8)$ codes, so they give at most $2^{10} = 1024$
vectors each by the linear programming bounds for binary codes.
The total number of possible vectors is $4600$ exactly, i.e., as
many as we started with. Therefore the numbers must be exact, and
in particular, we can normalize the code $\mathcal{D}$
corresponding to Case~\twoshref{III}, by the uniqueness of the $(3,6,22)$
Steiner system (Satz~4 in \cite{W}). We need to show that the
vectors of Case~\twoshref{IV} and \twoshref{V} are determined (up to isometries fixing
$V_0$, $\sqrt{2}D_{24}$, and the code $\mathcal{D}$) from this.
Let $W_0 = \left(3e_1 + e_2 - \sum_{i=3}^{24} e_i
\right)/\sqrt{8}$ be a vector from Case~\twoshref{IV}, which we can assume
after applying isometries as before. Let $V$ be a vector from
Case~\twoshref{III} with $\pm 2$'s at locations in the codeword $d \in
\mathcal{D}$, and suppose there are $r$ minus signs and $6-r$ plus
signs. Then
$$
\langle W_0, V \rangle  = \frac{1}{8} (6+2 + 2r -2(6-r)) =
\frac{-4+ 4r}{8},
$$
which forces $r$ to be odd. Now we get $2^6/2$ vectors for each
codeword, for $77$ codewords. Again exactness shows us that all
the vectors of Case~\twoshref{III} are uniquely determined.

Next we would like to show that vectors of Cases~\twoshref{I}, \twoshref{II} and~\twoshref{III}
span the lattice $L$, or in particular, the remaining generators
from Cases~\twoshref{IV} and~\twoshref{V}. It suffices to deal with Case~\twoshref{IV} since
clearly Case~\twoshref{V} is obtained by subtracting vectors of Case~\twoshref{IV} from
$4(e_1+e_2)/\sqrt{8} = V_0$.

For Case~\twoshref{IV}, we employ the same technique used in Case~\twoshref{IV} for the
$(22,891,1/4)$ code. It amounts to showing that the linear span of
the $77$ codewords $11d$ with $d \in \mathcal{D}$ contains exactly
$1024$ vectors of the form $00c$ with $c \in \F_2^{22}$, which is
easily checked on a computer.
\end{proof}

\section*{Acknowledgements}

We thank Eiichi Bannai for pointing out the reference \cite{Cu}
and the anonymous referee for helpful feedback.

\end{document}